\documentclass[12pt]{amsart}
\usepackage{amsmath,amssymb,amsfonts,amscd}
\usepackage[margin=1in]{geometry}
\usepackage{stmaryrd}
\usepackage{color}
\usepackage{scalerel}

\def\bal{\begin{aligned}}
\def\eal{\end{aligned}}
\def\be{\begin{equation}}
\def\ee{\end{equation}}
\def\bcs{\begin{cases}}
\def\ecs{\end{cases}}

\newtheorem{theorem}{Theorem}
\newtheorem{lemma}[theorem]{Lemma}
\newtheorem{proposition}[theorem]{Proposition}
\newtheorem{remark}[theorem]{Remark}
\newtheorem{definition}[theorem]{Definition}

\def\C{{\mathbb C}}
\def\Z{{\mathbb Z}}

\def\Q{{\mathbb Q}}
\def\R{{\mathbb R}}
\def\P{{\mathbb P}}

\def\lb{\llbracket}
\def\rb{\rrbracket}
\def\sH{\mathcal{H}}
\def\G{{\mathbf \Gamma}}

\title{Frobenius constants for families of elliptic curves}
\author{Bidisha Roy, Masha Vlasenko}

\address[]{Institute of Mathematics of the Polish Academy of Sciences \\ Śniadeckich 8, 00-656 Warsaw}
\email{m.vlasenko@impan.pl, broy@impan.pl}
\thanks{
Work of Masha Vlasenko was supported by the National Science Centre of Poland (NCN), grant UMO-2020/39/B/ST1/00940.}

\begin{document}
\date{\today}
\maketitle

\section{Introduction}

In this paper we will consider complex analytic functions which arise in the following way. Take a family of elliptic curves with one parameter, e.g. the Legendre family $y^2 = x (x-1)(x-t)$, $t \ne 0,1$. For this family there is a distinguished period integral which is analytic near $t=1$, namely
\[
\delta(t) = \int_t^1 \frac{dx}{\sqrt{x(x-1)(x-t)}}.
\]
We now choose a path from $t=1$ to another singular point $t=0$, consider the analytic continuation of $\delta(t)$ along this path and take its Mellin transform. Up to a simple factor of $s^2$, this yields the function of our interest: 
\be\label{kappa-Legendre}
\kappa(s) = s^2 \int_0^1 t^{s-1} \delta(t) dt = \sum_{n \ge 0} \kappa_n s^n.
\ee
The coefficients $\kappa_0,\kappa_1,\ldots$ of its power series expansion are called \emph{Frobenius constants}. Their study was initiated in~\cite{GZ1}, where the authors consider families of algebraic surfaces arising in the context of mirror symmetry. The general definition of Frobenius constants and their basic properties can be found in~\cite{BV}. These numbers are defined in terms of the differential equations satisfied by period integrals, so-called Picard--Fuchs differential equations of families of algebraic varieties. Period integrals of the Legendre family are annihilated by the hypergeometric differential operator
\be\label{hg-de}
t(1-t) \frac{d^2}{dt^2} + (1-2t) \frac{d}{dt} - \frac14.
\ee
Frobenius constants describe the monodromy of Frobenius solutions determined near one singularity along a path connecting it to another singularity, going around it and coming back,~\cite[Definition 22]{BV}.   Formula~\eqref{kappa-Legendre} is a special case of the main result of \emph{loc.cit.,}  \\cite[Theorem 30]{BV}, which connects Frobenius constants and expansion coefficients of generalized gamma functions.  We will recall the definition of generalized gamma functions in Section \ref{sec:proofs}.   By~\cite[Corollary 31]{BV} Frobenius constants of families of algebraic varieties are periods in the sense of Kontsevich--Zagier,~\cite{KZ-Periods}. Frobenius constants are related to periods of limiting Hodge structures, see~\cite{BV, Kerr}, but their true nature is yet   mysterious. For the moment we only know how to compute Frobenius constants for hypergeometric differential equations, see~\cite[Prop. 26]{BV}. In other cases one   may hope to evaluate a few constants $\kappa_n$ for small $n$. We will not assume that the reader is familiar with \emph{loc.cit.}. In the following sections we will view Frobenius constants simply as expansion coefficients of generalized gamma functions. Their relation to monodromy is only mentioned in this Introduction for motivational purpose. 

Our paper deals with Frobenius constants for \emph{stable} families of elliptic curves. A family is called stable if its singular fibres have only double points. Such families are known to have at least four singular fibres. Stable families with exactly four singular points were classified by Beauville,~\cite{Beau}. The respective Picard--Fuchs differential equations were found in~\cite{Ver}, they are of the second order with four singular points. In  Table~\ref{kappa-data} we will demonstrate results of numerical experiments in which some Frobenius constants $\kappa_n$  coincide to many decimal digits with $\Q$-linear combinations of products of zeta values of total weight $n$. Our data resembles the  conjectural evaluation of Frobenius constants for the Ap\'ery family of K3 surfaces in~\cite{GZ1}. In this paper we wanted to capture the simplest possible situation in which this curious phenomenon seems to take place. 

Families of elliptic curves admit a natural parametrization by modular functions. Their period integrals then correspond to modular forms of weight~1. In Sections~\ref{sec:results} and~\ref{sec:proofs} we will use this modular parametrization to express Frobenius constants as iterated integrals of modular forms in two different ways. The  first expression  given in Section~\ref{sec:results} is rather straightforward, it  uses the functional equation satisfied by generalized gamma functions. With this we will be able to evaluate Frobenius constants numerically. In Section~\ref{sec:proofs} we will use a more subtle, homological expression for generalized gamma functions, which allows us to represent  Frobenius constants $\kappa_n$ as regularized iterated integrals and to evaluate them for small $n$. The main results of our paper are Theorems~\ref{kappa-1-is-0-thm} and~\ref{kappa2-thm}.

In the end of this introduction let us compute Frobenius constants for the Legendre family. This turns out to be a simple task. We first change the order of integration in the integral expression in~\eqref{kappa-Legendre} and then substitute $t=xu$:
\[\bal
\int_0^1 \int_{t}^1 \frac{t^{s-1} dx dt}{\sqrt{x(1-x)(x-t)}} &= \int_0^1 x^{-1/2}(1-x)^{-1/2}\int_{0}^x t^{s-1} (x-t)^{-1/2} dt dx \\
&= \int_0^1 x^{-1/2}(1-x)^{-1/2} x^{s-1/2}\int_{0}^1 u^{s-1} (1-u)^{-1/2} du dx\\
& = \int_0^1 x^{s-1}(1-x)^{-1/2} dx \int_{0}^1 u^{s-1} (1-u)^{-1/2} du \\
&= \left( \frac{\Gamma(s)\Gamma(1/2)}{\Gamma(s+1/2)}\right)^2 = \pi\,\frac{\Gamma(s)^2}{\Gamma(s+\frac12)^2}.
\eal\]
This yields
\[
\kappa(s) = \pi \frac{s^2 \, \Gamma(s)^2}{\Gamma(s+\frac12)^2} = 16^s \frac{\Gamma(1+s)^4}{\Gamma(1+2s)^2} = \exp \left( 4 \log(2) s + 2 \sum_{k = 2}^\infty \frac{\zeta(k)}{k} (2 - 2^k)(-s)^k\right),
\]
where we first used the Legendre duplication formula for the classical gamma function and then the known power series expansion of the logarithm of the gamma function at $1$. Therefore we have
\be\label{hg-fc}
\kappa_0 =1, \; \kappa_1 = 4 \log(2), \; \ldots 
\ee
As we mentioned above, these constants are related to the monodromy of Frobenius solutions. The classical Frobenius method  gives the following standard solutions of the hypergeometric differential operator~\eqref{hg-de} near $t=0$:
\[\bal
&\phi_0(t) = \sum_{n=0}^\infty \frac{(\tfrac12)_n^2}{n!^2} t^n , \quad \phi_1(t) = \, \log(t) \phi_0(t) + \sum_{n=1}^\infty \frac{(\tfrac12)_n^2}{n!^2} \left(\sum_{k=1}^n \frac1{k(k-\tfrac12)}\right)  t^n .\\
\eal\]
The value $\kappa_1 = 4 \log (2)$ in~\eqref{hg-fc} means that the linear combination $\phi_1(t) - 4 \log(2) \phi_0(t)$ is analytic at $t=1$. This linear combination actually equals $-i \delta(t)$, see~\cite[(9)]{Vl-HS}. There are also higher Frobenius functions $\phi_2(t),\phi_3(t),\ldots$ which come as a natural extension of the Frobenius method, see~\cite{GZ1} and~\cite[\S 3]{BV}. For every $n$ the constant $\kappa_n$ is determined by the  property that the linear combination $\phi_n(t)-\kappa_n \phi_0(t)$ extends analytically through $t=1$. 

\smallskip

{\bf Acknowledgement.} We thank Francis Brown for explaining to us the concept of regularization which is used in this paper. We are also grateful to both our referees for a number of useful remarks and suggestions, which helped us to improve the exposition.

\section{Frobenius constants of Beauville's families}\label{sec:results}
The objective of this paper is the computation of Frobenius constants for the   six second order differential operators of the shape 
\be\label{L}
L = \theta^2  - t \bigl(A\,  \theta(\theta + 1) + \lambda \bigr) + B t^2 \,  (\theta+1)^2,
\ee
where $\theta=t \frac{d}{dt}$ and $A,B,\lambda$ are certain integer numbers whose values are listed in Table~\ref{Z-data}. Denote by $\phi_0(t)=\sum_{m \ge 0} u_m t^m$ the unique power series solution to $L$ normalized so that $u_0=1$. Observe that $L$ has four singular points whenever 
\be\label{distinct-roots}
B \ne 0 \text{ and } A^2 \ne 4B.
\ee
Under assumption~\eqref{distinct-roots}, the only known values $(A,B,\lambda)$ for which $\phi_0(t) \in \Z\lb t \rb$ are essentially those listed in Table~\ref{Z-data} and case {\bf B} with $(A,B,\lambda)=(9,27,3)$, see~\cite{Z-Apery}. Operators with parameters $(-A,B,-\lambda)$ and $(A,B,\lambda)$  are related via the change of variable $t \mapsto -t$. Adding such triples $(-A,B,-\lambda)$ for the mentioned cases {\bf A}-{\bf F} will give the full list of known operators $L$ with four singular points and $\phi_0\in \Z\lb t \rb$.   Each of these operators is a Picard-Fuchs operator associated to one of the stable families of elliptic curves considered by Beauville.  The correspondence between them and Beauville's classification can be found in~\cite[\S 7]{Z-Apery} and references therein. We will not need the equations of those families. What we will use is the modular parametrization, which we will also borrow from~\cite{Z-Apery}. Namely, for each case there exists a modular function $t(z)$ such that $f(z)=\phi_0(t(z))$ is a modular form of weight~1; here $z$ is the variable in the upper halfplane. In all cases except {\bf B} the value $c=t(0)$ is a singularity of the differential operator $L$. For brevity, we have excluded case {\bf B} from our considerations.

The singularities of $L$ are located at $t=0,\infty$ and the roots of $1-A t + B t^2=0$. One can check that each of them is a regular singular point with  maximally unipotent local monodromy transformation. The local exponents at $0$ and the roots of $1-A t + B t^2=0$ are equal to $0$, the local exponent at $\infty$ equals $1$. In particular, near each singularity $t=c$ there is a unique up to constant multiple solution $\delta^{(c)}(t)$ to $L$ which is analytic at $c$, and every other solution adds a constant multiple of $\delta^{(c)}(t)$ when going around this singular point. For example, one can take  the above-mentioned power series solution $\phi_0(t)=$   $\delta^{(0)}(t)$. 

\begin{table}
\caption{Modular parametrization of some Picard--Fuchs differential operators of Beauville's families of elliptic curves, \cite{Z-Apery}. In this table a standard notation for eta-function ratios is used, e.g. $\frac{a^kb^\ell}{c^m d^n}$ means $\frac{\eta(az)^k \eta(bz)^\ell}{\eta(cz)^m \eta(dz)^n}$ where $\eta( \cdot)$ is the Dedekind eta function.}
\begin{tabular}{c|ccc|cccc|c}
Zagier's &&&& Modular group& Hauptmodul & & &  $c=t(0)$ \\
label &$A$& $B$ & $\lambda$&$G$&$t(z)$&$f(z)$ & $\chi$&\\
\hline
&&&&&&&\\
\bf{A}&7&-8&2 & $\Gamma_0(6)$ & $\frac{1^36^9}{2^33^9}$ & $\frac{2^1 3^6}{1^26^3}$  & $\left(\frac{-3}{\cdot}\right)$ & $1/8$\\
&&&&&&&\\
\bf{C}&10&9&3& $\Gamma_0(6)$ & $\frac{1^46^8}{2^8 3^4}$ & $\frac{2^63^1}{1^36^2}$  & $\left(\frac{-3}{\cdot}\right)$ &1/9  \\
&&&&&&&\\
\bf{F}&17&72&6 & $h^{-1}\Gamma_0(6) h$ & $\frac{1^5 3^1 4^5 6^2 12^1}{2^{14}}$ &  $\frac{2^{15} 3^2 12^2}{1^64^66^5}$ & $\sim^h \left(\frac{-3}{\cdot}\right)$ & 1/9 \\
&&&&$\quad \begin{small} h=\begin{pmatrix}1&1/2\\0&1\end{pmatrix}\end{small}$&&& on $\Gamma_0(6)$\\
\hline
&&&&&&&\\
\bf{E}&12&32&4 & $  \Gamma_0(8)$ & $\frac{1^4 4^2 8^4}{2^{10}}$ & $\frac{2^{10}}{1^4 4^4}$ & $\left(\frac{-4}{\cdot}\right)$ &1/8 \\
&&&&&&&\\
\bf{G}&0&-16&0 & $  \Gamma_0(8)$  & $\frac{2^4 8^8}{4^{12}}$ & $\frac{4^{10}}{2^4 8^4}$ & $\left(\frac{-4}{\cdot}\right)$ & 1/4 \\
&&&&&&&\\
\hline
&&&&&&&\\
\bf{D}&11&-1&3 &$ \Gamma_1(5)$ & $ q \prod_{n=1}^\infty (1-q^n) ^{ 5 (\frac{n}{5})}$ & $\left(t^{-1} \frac{5^5}{1^1}\right)^{1/2}$  & $1$ & $-\frac{11}2+\frac52\sqrt{5}$ \\
&&&&&&&\\
\end{tabular}
\label{Z-data}
\end{table}

\begin{proposition}\label{gamma-properties} In all cases listed in Table~\ref{Z-data} the function defined for $\mathrm{Re}(s)>0$ by the integral
\be\label{gen-gamma}
\G(s) =  \int_{0}^{i \infty} t(z)^{s-1} z f(z) d t(z)
\ee
 satisfies the functional equation
\be\label{func-eq}
s^2 \G(s) = (As(s+1)+\lambda)\G(s+1) - B (s+1)^2 \G(s+2).
\ee
\end{proposition}

%Note that, similarly to the classical gamma function, one can extend $\G(s)$ to a meromorphic function in the complex plane using the the above functional equation. In the next section we will give a different expression for $\G(s)$ which will show that it is in fact an entire function. For this property the factor of $(e^{2 \pi i s}-1)^2$ in front of the integral is essential.

\begin{proof}
 Consider the singularity $c=t(0)$ of the respective differential operator $L$. It is a property of the modular parametrization in Table~\ref{Z-data} that under the substitution $t=t(z)$ the space of solutions of $L$ consists of $\C$-linear combinations of $f(z)$ and $z f(z)$. One can check that in each case $z f(z)$ has a finite limit when $z \to 0$. It follows that 
$z f(z) = \delta^{(c)}(t(z))$ for a solution $\delta^{(c)}(t)$ of $L$ which is analytic near $t=c$. Indeed, as $t=c$ is a point of maximally unipotent monodromy for $L$ with local exponent~$0$ then non-analytic at $t=c$ solutions grow infinitely as $t \to c$. We then have
\begin{equation}\label{Gamma-fun}
\G(s) =  \int_c^0 t^{s-1} \delta^{(c)}(t) dt.
\end{equation}
Next we consider  the identity
$ \displaystyle \int _{0}^{c} t^{s-1} L(\delta^{(c)}(t)) dt =0$.  Integrating by parts we obtain
\begin{align*} 
&s^2 \int _{0}^{c} t^{s-1} \delta^{(c)}(t) dt - A( s(s+1) + \lambda)\int _{0}^{c} t^{s} \delta^{(c)}(t) dt +B(s+1)^2 \int _{0}^{c} t^{s+1} \delta^{(c)}(t) dt\\
& = s c^s \delta ^{(c)}(c) ( 1 -Ac+ Bc^2) - c^{s+1} \frac{d}{dt} \delta^{(c)}(t) \mid_{t=c} (1- Ac +Bc^2).
\end{align*}
 Since $c$ is a root of $ 1 - At + Bt^2$, the  right-hand side vanishes. Thus  we obtain the desired functional equation of $ \G(s)$, namely \eqref{func-eq}.  
\end{proof}

Note that, similarly to the classical gamma function, one can extend $\G(s)$ to a meromorphic function in the complex plane using the above functional equation~\eqref{func-eq}.
The \emph{Frobenius constants} that we want to compute are given by the expansion coefficients at $s=0$ of the function
\be\label{kappa}
\kappa(s) =   2 \pi i \, s^2 \, \G(s) = \sum_{n \ge 0} \kappa_n s^n. 
\ee 
Note that $t(z)$ and $f(z)$ take real values on the imaginary half-axis, which can be seen from their eta-product presentations in Table~\ref{Z-data}. Therefore it follows from~\eqref{gen-gamma} and~\eqref{func-eq} that $\kappa(s)$ is real-valued for $s \in \R$.
 It follows from~\eqref{func-eq} that
\be\label{kappa-via-func-eq}
\bal
\kappa(s) = 2 \pi i & (As(s+1)+\lambda) \int_0^{i \infty} t(z)^{s} z f(z) t'(z) dz \\
& - 2 \pi i B (s+1)^2 \int_0^{i \infty} t(z)^{s+1} z f(z) t'(z) dz,
\eal\ee
where the expression in the right-hand side is well-defined for $\mathrm{Re}(s)>-1$. In particular, we see that $\kappa(s)$ is analytic at $s=0$.
  
 We differentiate~\eqref{kappa-via-func-eq} in $s$   and subsequently evaluate  it at $s=0$, which  yields an expression for $\kappa_n$ involving integrals
\[
\int_0^{i \infty} \log^m( t(z)) z  f(z) t'(z) dz, \quad \int_0^{i \infty}  t(z)\log^m( t(z)) z   f(z) t'(z) dz
\]  
with $m \le n$.  We computed   these integrals numerically in PARI/GP with a precision of about  75 digits after the comma.  Remarkably, we were often able to identify  $\kappa_n$ as $\Q$-linear combinations of products of zeta values of total weight $n$. The respective conjectural identifications of $\kappa_n$ are listed in Table~\ref{kappa-data}.  In case {\bf G} this fact can be actually proved for all $n$, see the argument in the end of this section. However in the cases marked by * in  Table \ref{kappa-values} such an identification  doesn't seem possible already for $n=4$, in which case one could expect rational multiples of $\pi^4$.  These Frobenius constants  may be periods of more complicated nature.

\begin{table} \label{kappa-values}
\caption{Conjectural  values (observed equal numerically up to 75 digits) of  some Frobenius constants $\kappa_n$. Entries * mean that we couldn't identify the respective constants as linear combinations of products of zeta values. The values of $\kappa_0,\kappa_1$ and $\kappa_2$ as well as other values for the case {\bf G} are proved in this paper, while the values in the upper-right corner remain conjectural.   }
\begin{tabular}{c|ccc|ccc}
case&$\kappa_0$&$\kappa_1$&$\kappa_2$&$\kappa_3$&$\kappa_4$&$\kappa_5$\\
\hline
\bf{A}&$1$&$0$&$-\frac{\pi^2}6$&$\zeta(3)$&$\frac{\pi^4}{80}$ & $-\frac7{16}\zeta(5)-\frac{5}{16}\pi^2\zeta(3)$ \\
&&&&&\\
\bf{C}&  $1$ &$0$ & $-\frac{\pi^2}6$ & $ \frac{2}{3} \zeta(3)$ & * & * \\
&&&&&\\
\bf{D}& $1$ & $0$ & $ -\frac{7}{30}\pi^2$ & $ 2 \zeta(3)$ & $ \frac{\pi^4}{180}$ & $\zeta(5) - \frac{\pi^2}2 \zeta(3)$ \\
&&&&&\\
\bf{E}& $ 1$ & $0$ & $ - \frac{\pi^2}{12}$ & $ - \frac{5}{4} \zeta(3)$ & * & *\\
&&&&&\\
\bf{F}& $1$ & $0$ & $0$ & $ - \frac{11}{3} \zeta(3)$ & * & *\\ \cline{5-7}
&&&\multicolumn {1}{ c }{}&&\\
\bf{G}& $1$ & $0$ &  \multicolumn {1}{ c }{ $- \frac{\pi^2}{12}$ } & $ \frac{1}{2} \zeta(3)$ & $ -\frac{\pi^4}{720}$ & $\frac38 \zeta(5) - \frac{\pi^2}{24}\zeta(3)$\\
\end{tabular}
\label{kappa-data}
\end{table}

In the following sections we will be able to confirm the values of $\kappa_0$, $\kappa_1$ and $\kappa_2$ using a homological interpretation of the function $\G(s)$. To state our main results we will need the following concept of \emph{regularization at infinity} for functions in the upper halfplane. We learned it from Francis Brown, see e.g.~\cite[\S4-5]{Brown-MMV}. In Section~\ref{sec:kappa2} we shall see how this definition is compatible with the classical theory of periods of modular forms.

\begin{definition}\label{regularized integral} We say that a function $F:\sH \to \C$ is \emph{regularizable} as $z \to \infty$ if there exists a polynomial $P \in \C[z]$ such that
\[
F(z) = P(z) + o( \mathrm{Im}(z)^{-N}), \quad \text{for all } N \ge 1
\] 
when  $\mathrm{Re}(z)$ is fixed and $\mathrm{Im}(z) \to \infty$. We define the regularized value of $F(z)$ at infinity as
\[
[F(z)]_{reg} := P(0).
\]
\end{definition}  

\begin{theorem}\label{kappa-1-is-0-thm} In cases {{\rm \textbf{A, C-G}}} given in Table~\ref{Z-data} the expansion coefficients of the function
\[
\kappa(s) = 2\pi i s^2 \int_0^\infty t(z)^{s-1} z f(z) d t(z) = \sum_{n=0}^\infty \kappa_n s^n
\]
are given by  $\kappa_0=1$, $\kappa_1=0$ and
\be\label{kappa-via-iterated-int}
\kappa_{2+m} = \frac{4 \, \pi^2}{b \, m!} \left[\int_{z_0/(b z_0+1)}^{z_0} \log^m (t(z)) g(z) dz\right]_{reg} 
\ee
for all $m \ge 0$, where $g(z) = \frac1{2 \pi i} \frac{t'(z)}{t(z)} f(z)$ is a modular form on $G$ of weight~3 and character $\chi$ and $b$ is the minimal positive integer such that $\begin{pmatrix}1 & 0 \\ b & 1\end{pmatrix} \in G$. 
\end{theorem}

For convenience let us list the values of $b$ in this theorem:  

\bigskip
\begin{center}
\begin{tabular}{c|cccccc} \label{values of 'b'}
case & {\bf A} & {\bf C} & {\bf D} & {\bf E} & {\bf F} &  {\bf G} \\
&&&&&&\\
\hline
&&&&&&\\
$b$ &6 & 6 & 5 & 8 & 12 & 8\\
\end{tabular}
\end{center}

\bigskip

Theorem~\ref{kappa-1-is-0-thm} will be proved in Section~\ref{sec:proofs}. In Section~\ref{sec:kappa2} we will evaluate the regularized integral for $m=0$ and confirm the values of $\kappa_2$: 

\begin{theorem}\label{kappa2-thm}   The  values of $\kappa_2$ given in Table~\ref{kappa-data}  are equal to  the  regularized
integral expression given  in Theorem \ref{kappa-1-is-0-thm} for $m =0$. 
\end{theorem}

At the end let us remark that all Frobenius constants in case {\bf G} can be evaluated using the fact that in this case $L=\theta^2-16 t^2(\theta+1)^2$ is  the change of variable $t \mapsto 16 t^2$ in the hypergeometric differential equation~\eqref{hg-de}. Therefore $\delta^{(c)}(t)$ in this case is a constant multiple of $\delta(16 t^2)$ and by simple manipulation with integrals we find that
\[
s^2 \int_0^1 t^{s-1} \delta(t) dt = 2 s^2 16^s \int_0^{1/4} t^{2s-1} \delta(16 t^2) dt = 16^s \kappa_{\text{\bf G}}(2s). 
\]
This function was computed at the end of the introductory section, and therefore we obtain
\begin{equation} \label{kappa values of G}
\kappa_{\text{\bf G}}(s) = \exp\left( 2 \sum_{k = 2}^\infty \frac{\zeta(k)}{k} (2^{1-k} - 1)(-s)^k \right) = 1 - \frac{\pi^2}{12} s^2 + \frac{\zeta(3)}{2} s^3 + \ldots
\end{equation}

\section{Generalized gamma functions}\label{sec:proofs}

For a differential operator $L$ on $\P^1(\C)$ with regular singularities, one can consider Mellin transforms of its solutions along closed loops, e.g.
\[
\G_{(\sigma, \phi)}(s) = \int_\sigma t^{s-1} \phi(t) dt
\]
The natural requirements here are that $\sigma$ is a closed loop avoiding the singularities of $L$, $\phi(t)$ is a solution to $L \phi=0$ having no monodromy along the loop $\sigma$. Moreover, we should require that $t^s$ is single-valued along $\sigma$. This latter restriction is equivalent to the fact that $\sigma$ is contractible in $\C^* = \C \setminus \{0\}$. Such a choice of $\sigma$ and $\phi$ may not always exist. More generally, one can attach gamma functions to homology classes with coefficients in the twisted local system of solutions of $L$, see~\cite[\S 1]{BV}. To describe such homology classes it is convenient to fix a regular point $p \in \C^*$ and consider the group $G=\pi_1(\P^1(\C) \setminus S ,p)$ of homotopy classes of loops based at $p$. Here $S$ is a finite set of points consisting of $0,\infty$ and the singularities of $L$. Let $V$ be the $\C$-vector space of solutions of $L$ near $p$. Then $V$ is  a representation $G$ by monodromy of solutions along loops. Consider the ring 
\[
R=\C[e^{\pm 2 \pi i s}]
\]
and fix some branch of $t^s$ near the point $t=p$. Then $R$ is a representation of $G$ via the monodromy of $t^s$. A 1-cycle $\xi$ for the homology of $G$ in the twisted representation $\tilde V = V \otimes_\C \C[e^{\pm 2 \pi i s}]$ is a finite sum
\[
\xi = \sum_j \sigma_j \otimes \phi_j \otimes e^{2 \pi i s n_j}
\]
where $\sigma_j \in G$, $\phi_j \in V$, $n_j \in \Z$ and   
\[
\partial \xi = \sum_{j}\left( \sigma_j\left( \phi_j \otimes e^{2 \pi i s n_j} \right) - \phi_j \otimes e^{2 \pi i s n_j} \right) = 0.
\]   
The respective generalized gamma function is defined as
\[
\G_\xi(s) = \sum_je^{2 \pi i s n_j} \int_{\sigma_j} t^{s-1} \phi_j(t) dt.
\]
This is an entire function of $s$, which depends only on the homology class of $\xi$ in $H_1(G,\tilde V)$ and satisfies a functional equation, \cite[Lemma~4 and Proposition~8]{BV}. For our differential operator~\eqref{L} the functional equation that generalized gamma functions satisfy reads as~\eqref{func-eq}.  

\begin{remark} Let us mention that generalized gamma functions corresponding to a differential operator $L$ form an $R$-module of finite rank, \cite[Proposition 5]{BV}. This can be seen in the following way. Observe that $\tilde V$ is an $R$-module and the action of $G$ is $R$-linear. Therefore the homology group $H_1(G,\tilde V)$ is an $R$-module as well. The module of gamma functions is the quotient of $H_1(G,\tilde V)$ by those homology classes whose gamma functions vanish.

Though we will not use this fact, let us mention that the $R$-module of gamma functions of our differential operator~\eqref{L} has rank at most~2. The rank drops to~1  when $A=\lambda=0$, in which case $L$  comes from a hypergeometric differential operator by a change of variables. It was shown in~\cite[\S 2]{BV} that for hypergeometric differential operators  the module of generalized gamma functions has rank~1. 
\end{remark}

We will now consider generalized gamma functions in a situation when $L$ is of order~2 and admits a modular parametrization. Namely, let $\sH$ be the upper halfplane and  $t: \sH \to \C^*$ be a Hauptmodul for some genus zero  congruence subgroup $G \subset SL_2(\Q)$. Let $S \subset \P^1(\C)$ be the set of cuspidal values of $t(z)$. Let $f(z)$ be a modular form $f(z)$ of weight one on $G$ with a character $\chi$. This function satisfies 
\[
f\left(\frac{az+b}{cz+d}\right) = \chi(\sigma) (c z + d) f(z)
\]
for any $\sigma= \begin{pmatrix}a & b \\ c & d\end{pmatrix} \in G$.  
 In this situation there exists a differential operator $L$ of order~2 on $\P^1(\C) \setminus S$ whose pull-back under $t(z)$ annihilates the 2-dimensional $\C$-vector space spanned by  $f(z)$ and $z f(z)$, ~\cite[\S 5.4]{1-2-3}. 

We choose the base point $p \in \P^1(\C) \setminus S$ and identify the space of solutions of $L$ near $p$ with the space of polynomials $V=\C z + \C$. Here $\sigma= \begin{pmatrix}a & b \\ c & d\end{pmatrix} \in G$ acts on a polynomial $Q(z) \in V$ by sending it into the polynomials $\chi(\sigma) (c z + d) Q\left(\frac{az+b}{cz+d}\right)$.

To describe $\tilde V$, we need to twist by the monodromy of $t^s$. We assume that $t(z)$ has no zeroes in $\sH$. Take any branch of $\log t(z)$ and note that it is single-valued in $\sH$. Then for $\sigma \in G$ we define 
\[
n(\sigma) = \frac1{2 \pi i} \left( \log(t(\sigma z))-\log t(z) \right) \in \Z.
\]
This is an integer which is independent of $z$ because $\sH$ is connected. One can easily see that $n(\sigma_1 \sigma_2)=n(\sigma_1)+n(\sigma_2)$. Then $\tilde V$ is the twist of $V$ by the multiplicative character $\gamma \mapsto e^{2 \pi i s n(\gamma)} \in \C[e^{\pm 2 \pi i s}]^\times$. We recall the notation $R=\C[e^{\pm 2 \pi i s}]$ and identify $\tilde V \cong R z + R$. The action of $\sigma= \begin{pmatrix}a & b \\ c & d\end{pmatrix} \in G$ on an element  $Q \in \tilde V$ is given by formula  
\[
(Q | \sigma )(z) = e^{2 \pi i s n(\sigma)} \chi(\sigma) (c z + d) Q\left(\frac{az+b}{cz+d}\right). 
\]
Our notation $Q|\sigma$ reflects the fact that this is a right action. 

A homological 1-cycle with coefficients in $\tilde V$ is a finite sum $\xi = \sum_j \sigma_j \otimes Q_j$ satisfying the condition
\[
\partial \xi = \sum_j \,\left( Q_j | \sigma_j - Q_j \right) = 0.  
\]  
The respective generalized gamma function is given by 
\be\label{modular-gamma-1}
\G_{\xi}(s) = \sum_j \int_{z_0}^{\sigma_j z_0} t(z)^{s}Q_j(z)g(z) dz,
\ee
 where
\be\label{g-wt-3-form}
g(z) = \frac1{2 \pi i} f(z) \frac{t'(z)}{t(z)}
\ee
is a modular form on $G$ of weight~3 and character $\chi$. Here $z_0 \in \sH$ is  an arbitrary point in $\sH$, it is a preimage of the basepoint $p=t(z_0)$. Note that the above integral~\eqref{modular-gamma-1} is well-defined for any $s \in \C$, and therefore it defines an entire function $\G_{\xi}(s)$.  The following properties are  expected but worth to check. 

\begin{lemma}\label{indep-of-z0} The function  $\G_\xi(s)$  defined by \eqref{modular-gamma-1} is independent of $z_0 \in \sH$. This function depends only on the class of $\xi$ in  the homology group $H_1(G,\tilde V)$.
\end{lemma}
\begin{proof}To prove independence of $z_0$ we transform the difference 
\[\bal
&\sum_j \left( \int_{z_0}^{\sigma_j z_0} t(z)^{s} Q_j(z) g(z) dz - \int_{z_1}^{\sigma_j z_1} t(z)^{s} Q_j(z) g(z) dz \right) \\
&= \sum_j \left( \int_{z_0}^{z_1} t(z)^{s} Q_j(z) g(z) dz - \int_{\sigma_j z_0}^{\sigma_j z_1} t(z)^{s} Q_j(z) g(z) dz \right)
\eal\]
and perform the change of variable $z=\sigma_j u$ in the second integral in the brackets. Since $g(z)$ is a modular form of weight~3 on $G$ with character $\chi$, then for any $\sigma= \begin{pmatrix}a & b \\ c & d\end{pmatrix} \in G$ and $Q \in \tilde V$ one has
\[
t(z)^{s} Q(z) g(z) dz |_{z = \sigma u} = t(u)^{s} (Q|\sigma)(u) g(u) du.
\]
With this observation, the above difference becomes
\[
= \int_{z_0}^{z_1} t(z)^{s}\sum_j \left(Q_j(z) - Q_j|\sigma_j(z)\right) g(z) dz = - \int_{z_0}^{z_1} t(z)^{s} (\partial \xi)(z) g(z) dz = 0.
\]
Here we used the fact that $\partial \xi = 0$. For the second statement we evaluate $\G_\xi(s)$ for a 1-boundary 
\[
\xi = \partial \; \left( (\sigma_1,\sigma_2) \otimes Q \right) = \sigma_2 \otimes (Q|\sigma_1) - \sigma_1 \sigma_2 \otimes Q + \sigma_1 \otimes Q.
\]
One has
\[\bal
\int_{z_0}^{\sigma_1 \sigma_2 z_0}  t(z)^{s} Q(z) g(z) dz &= \int_{z_0}^{\sigma_1 z_0}  t(z)^{s} Q(z) g(z) dz + \int_{\sigma_1 z_0}^{\sigma_1 \sigma_2 z_0}  t(z)^{s} Q(z) g(z) dz \\
&= \int_{z_0}^{\sigma_1 z_0}  t(z)^{s} Q(z) g(z) dz + \int_{z_0}^{\sigma_2 z_0}  t(z)^{s} (Q|\sigma_1)(z) g(z) dz, \\
\eal\]
which precisely means that $\G_\xi(s) = 0$. 
\end{proof}

We now turn to the proof of   Theorem~\ref{kappa-1-is-0-thm}. From now on we assume that $L$ is of the shape~\eqref{L} with $(A,B,\lambda)$, $t(z)$ and $f(z)$ listed in Table~\ref{Z-data}. Note that in all our cases $t(i \infty)=0$. We consider the singularity of $L$ given by the cuspidal value $c=t(0)$ and the direct path from $0$ to $c$ along the real line. It was noticed in~\cite[\S 2]{BV} that if we restrict gamma functions to a neighbourhood of this path in $\P^1(\C) \setminus S$, there will be essentially a unique such function to consider. This is  the function $\G(s)$ that was introduced in formula~ \eqref{gen-gamma}. More precisely, we will show that 
\be\label{two-gammas}
\G_{\xi}(s) =  \frac{(e^{2 \pi i s}-1)^2}{2 \pi i} \; \G(s), 
\ee
where $\xi$ is a homology 1-cycle whose class generates  the quotient 
\be\label{homology-quotient}
H_1(\langle\sigma_0,\sigma_c\rangle, \tilde V) / H_1(\langle \sigma_c \rangle, \tilde V)
\ee
as a $\C[e^{\pm 2 \pi i s}]$-module.  Here by $\langle .. \rangle$ we denote the subgroup of $G$ generated by the respective elements and  $\sigma_0,\sigma_c$ are elements of $G$ which in $t$-plane correspond to the loops around $0$ and $c$ respectively.   We take for $\sigma_0$ and $\sigma_c$  the  generators of  stabilizers of the respective cusps $z=i\infty$ and $z=0$ in $G$. In all our cases these generators are of the form
\[
\sigma_0 = \begin{pmatrix}1 & 1 \\ 0 & 1\end{pmatrix}, \qquad \sigma_c =  \begin{pmatrix}1 & 0 \\ b & 1\end{pmatrix},
\]
 where $b$ is precisely the integer from Theorem~\ref{kappa-1-is-0-thm}. It is easy to see that gamma functions corresponding to homology classes  $\mu \in H_1(\langle \sigma_c \rangle, \tilde V)$ vanish. Indeed, such gamma functions ${\bf \Gamma_{\mu}}(s)$ can be written as $\C[e^{\pm 2 \pi i s}]$-linear combinations of integrals of the form $\oint t^{s-1} \phi(t) dt $ where the integral is over a small circle around $c$ and $\phi(t)$ is a solution of $L$ analytic at $t=c$. Since $t^s$ is also analytic at $t=c$, we have $\G_\mu(s) \equiv 0$.

\begin{lemma}\label{explicit second order polynomial} The quotient~\eqref{homology-quotient} is a 
$\C[e^{\pm 2 \pi i s}]$-module of rank~1 generated by the class of
\[
\xi = \sigma_0 \otimes Q_0 + \sigma_c \otimes Q_c. 
\]
with $Q_0 = (1-e^{2 \pi i s} )z + e^{2 \pi i s}$, $Q_c = \frac1{b} (1-e^{2 \pi i s})^2$.
\end{lemma}

\begin{proof}
 Note that $\chi(\sigma_0)=\chi(\sigma_c)=1$ for the character $\chi$ in Table~\ref{Z-data}. Therefore, for proving this lemma, we assume $\chi = 1$.
Let $ Q_0(z) = \kappa z+  \mu$ and $ Q_c(z) = uz+v$, for some $ \kappa, \mu, u,v \in  \mathbb{C}[ e^{ \pm 2 \pi i s}]$ such that $ \sigma_0 \otimes Q_0 + \sigma_ c \otimes Q_c$  is an arbitary element of $ H_1 ( \langle \sigma_0, \sigma_c \rangle, \tilde{V})$.  We can assume $u=0$ because $ \partial(  \sigma_c \otimes z) = z-z=0$ and thus $  \sigma_c \otimes z \in H_1 ( \langle \sigma_c \rangle , \tilde{V})$. Therefore, we need $\kappa,t, v \in \mathbb{Q}[e^{ \pm 2 \pi i s}]$ such that
\begin{align*}
&  0 = \partial(\sigma_0 \otimes Q_0+ \sigma_c \otimes Q_c) =e^{2 \pi i s} (\kappa(z+1) +\mu) - (\kappa z+\mu) + (b z+1) v-v  \\
&\; =z ( e^{2 \pi i s} \kappa -\kappa +b v) +( e^{2 \pi is}(\kappa+\mu)-\mu).
\end{align*} It gives  $\displaystyle  e^{2 \pi i s} \kappa = \mu ( 1 - e^{ 2 \pi i s})$ and $ \displaystyle b v = \kappa ( 1 - e^{ 2 \pi is})$. Here $\mu \in \mathbb{C}[e^{\pm 2 \pi i s}]$ can be chosen arbitrarily, which shows that the quotient module is of rank one. Taking the unit value $\mu=e^{2 \pi i s}$ we obtain the polynomials $Q_0$ and $Q_c$ given in the statement of the lemma.
\end{proof}

 The gamma function   corresponding to the homological cycle $\xi$ in Lemma~\ref{explicit second order polynomial} is given by
\be\label{homological-gamma}
\G_{\xi}(s) = \int_{z_0}^{z_0+1} t(z)^{s} Q_0(z)g(z)dz + \int_{z_0}^{z_0/(b z_0 + 1)} t(z)^{s} Q_c(z)g(z)dz.   
\ee
To prove  relation \eqref{two-gammas}  we break each of the integrals into two parts as follows:  
\begin{align*}
\G_{\xi}(s) & =  \int_{z_0}^{i\infty}  Q_0(z) t(z)^s g (z) d z+ \int_{i\infty}^{z_0+1} Q_0(z) t(z)^s g(z) dz \\
& +  \int_{z_0}^{0} Q_c(z) t(z)^{s} g(z)dz+  \int_{0}^{z_0/(b z_0 + 1)} Q_c(z)t(z)^{s} g(z)dz\\
& = \int^{z_0}_{i\infty} \left( e^{ 2 \pi i s } Q_0(z+1) - Q_0(z) \right) t(z)^s g(z) dz \\
& \qquad +\int_0^{z_0} \left( (bz+1)Q_c(\frac{z}{bz+1}) - Q_c(z)\right)t(z)^s g(z) dz\\
& = -(1 - e^{ 2 \pi i s})^2  \int^{z_0}_{i\infty} z t(z)^s  g(z) dz \\
& \qquad +  \frac{1}{b} ( 1 - e^{2 \pi i s })^2 \int^{z_0}_{0}  ((bz +1) -1)t(z)^{s} g(z)dz\\
 &=  (1 - e^{ 2 \pi i s})^2  \int^{i\infty}_0 z t(z) ^s g(z) dz = \frac{(1-e^{2 \pi i s})^2}{2 \pi i} \, \G(s).
\end{align*}

To prove our theorem we need to compute the Taylor expansion at $s=0$  of the function 
\[
\kappa(s) = \left(\frac{2 \pi i s}{1-e^{2 \pi i s}}\right)^2 \G_{\xi}(s) = 2 \pi i s^2 \G(s).
\]
Note that each of the two integrals in the right-hand side of~\eqref{homological-gamma} may depend on $z_0$ while their sum is independent of $z_0$. Consider the Taylor expansion of the first integral at $s=0$. It turns out that the coefficients of this expansion have polynomial asymptotics when $z_0 \to \infty$.

\begin{lemma}\label{T-part-regularized} In each of our  six cases one has 
\[
 \left[\int_{z_0}^{z_0+1} t(z)^s g(z)  Q_0(z) dz \right]_{reg} =  \left( \frac{ 1- e^{2 \pi i s}}{2 \pi i s}\right)^2.
 \]
Here the integral is first expanded as a formal power series in $s$ and then each term is regularized as a function of $z_0$. 
\end{lemma}
\begin{proof}
 We denote $q=e^{2 \pi i z}$. After writing $ \log t(z) = 2 \pi i z + \sum_{  k \geq 1} c_k q^k$ and $ g(z) = \sum_{n \geq 0} a_n q^n$, we have 
\[
t(z)^s \cdot g(z) = \sum_{m \geq 0} \frac{s^m}{m!} \left( 2 \pi i z + \sum_{k \geq 1} c_k q^k\right)^m  \left(\sum_{n \geq 0} a_nq^n \right).
\]
The term near each $s^m$ a polynomial in $z$ whose coefficients are convergent $q$-series. Observe that  for any fixed value of $\mathrm{Re}(z_0)$ one has
\[
\int_{z_0}^{z_0+1} z^k (\sum_{n \ge 0} b_n q^n ) dz = b_0 \int_{z_0}^{z_0+1} z^k dz +  o( \mathrm{Im}(z_0)^{-N}), \quad  \text{for all } N \ge 1
\]
 when $\mathrm{Im}(z_0) \to \infty$. Therefore
\[
\left[\int_{z_0}^{z_0+1} z^k (\sum_{n \ge 0} b_n q^n ) dz \right]_{reg} = b_0 \int_0^1 z^k dz = b_0/(k+1).
\]
Hence, as equality of formal power series in $s$, we have 
\[
 \left[ \int_{z_0}^{z_0+1} t(z)^s g(z) Q_0(z) dz \right]_{reg} = a_0 \int_0^1 e^{2 \pi i zs} Q_0(z)dz = a_0 \left( \frac{ 1- e^{2 \pi i s}}{2 \pi i s}\right)^2. 
\]
Observe that $g(z)$ in~\eqref{g-wt-3-form} has $a_0=1$ for all cases under consideration, which concludes the proof of the lemma.
\end{proof}

\begin{proof}[Proof of Theorem~\ref{kappa-1-is-0-thm}.]
By Lemma~\ref{indep-of-z0} the sum of integrals in~\eqref{homological-gamma} is independent of $z_0$. Using Lemma~\ref{T-part-regularized} for the first integral, we conclude that each coefficient in the expansion of the second integral considered as a power series in $s$ is regularizable. With the expression for $Q_c$ from Lemma~\ref{explicit second order polynomial} we conclude that
\[\bal
\kappa(s) = \left( \frac{2 \pi i s}{ 1- e^{2 \pi i s}}\right)^2 \G_{\xi}(s) &= 1 - \frac{4 \pi^2 s^2}{b} \left[\int_{z_0}^{z_0/(b z_0+1)} t(z)^s g(z)  dz \right]_{reg} \\
&= 1 + \frac{4 \pi^2}{b} \sum_{m \ge 0} \frac{s^{m+2}}{m!} \left[\int^{z_0}_{z_0/(b z_0+1)} \log^m(t(z))  g(z)  dz \right]_{reg}. \\
\eal\]
\end{proof}

\section{Evaluating regularized integrals}\label{sec:kappa2}

In this section we prove Theorem~\ref{kappa2-thm}. What we need for this is to evaluate the regularized integral
\be\label{reg-int-to-do}
\left[\int_{z_0}^{z_0/(b z_0 + 1)} g(z) dz\right]_{reg}.
\ee
  Here, corresponding to each $t(z)$ and $f(z)$ as in Table \ref{Z-data}, we consider $g(z) = \frac1{2 \pi i} \frac{t'(z)}{t(z)} f(z)$ which becomes a weight~3 modular form on a modular group $G \subset SL_2(\Z)$ with character $\chi$ and $b$ is the minimal positive integer such that $\begin{pmatrix}1 & 0 \\ b & 1\end{pmatrix} \in G$. We will see later that, for proving Theorem~\ref{kappa2-thm}, it is enough to calculate the values of \eqref{reg-int-to-do} only for two cases, namely \textbf{A} and \textbf{D}. Therefore by possibly restricting $G$ to a smaller subgroup we may assume that $\chi=1$, which will simplify the considerations in what follows. 

We first review some properties of  integrals~\eqref{reg-int-to-do} and their relation to periods of modular forms. Let $g \in M_k(G)$ be a modular form of weight~$k$. For any $z_0$ in the upper halfplane $\sH$, consider the map $\psi_{z_0}: G \rightarrow V_{k-2}$ from $G$ into the set $V_{k-2} \subset \mathbb{C}[X]$ of polynomials of degree $\le k-2$ defined as 
\be\label{cocycle-z-0}
\psi_{z_0}(\gamma) = \int_{z_0} ^{\gamma z_0} (X-z) ^{k-2} g(z) dz.
\ee
The space $V_{k-2}$ is a representation of $G$ under the usual right action in weight $2-k$ which is given for $Q(X) \in V_{k-2}$ and $\gamma=\begin{pmatrix}a & b \\ c & d\end{pmatrix} \in G$ by $(Q|\gamma)(X) = (c X + d)^{k-2} Q(\frac{aX+b}{cX+d}) \in V_{k-2}$. We denote $\gamma Q = Q | \gamma^{-1}$ to turn it into the left action. The map~\eqref{cocycle-z-0} then satisfies the relation $ \psi_{z_0}(\gamma_1 \gamma_2) = \gamma_1\psi_{z_0}(\gamma_2) + \psi_{z_0}(\gamma_1)$, which precisely means that $\psi_{z_0}$ is a $1$-cocyle for $G$ with coefficients in the representation $V_{k-2}$. Note that
\be\label{cocycle-dependence-on-z0}\bal
\psi_{z_0}(\gamma) - \psi_{z_1}(\gamma) &= \int_{z_0}^{z_1} (X-z)^{k-2} g(z) dz  - \int_{\gamma z_0}^{\gamma z_1} (X-z)^{k-2} g(z) dz \\
&= (1- \gamma) \int_{z_0}^{z_1} (X-z)^{k-2} g(z) dz
\eal\ee
is a coboundary, and therefore we conclude that the class of $\psi_{z_0}$ in $H^1(G,V_{k-2})$ is independent of $z_0$. We denote this class by $[g] \in H^1(G,V_{k-2})$, it depends only on the modular form $g$.

\begin{proposition}\label{regularization-prop} When $\mathrm{Re}(z_0)$ is fixed and $\mathrm{Im}(z_0)\to \infty$ we have 
\be\label{cocycle-asymptotics}
\psi_{z_0}(\gamma) = \sum_{i=0}^{k-1} P_{i,\gamma}(X)z_0^i + o(\mathrm{Im}(z_0)^{-N}), \; \text{for all } \; N \ge 1
\ee
with some  polynomials $P_{i,\gamma} \in V_{k-2}$ that are independent of $z_0$. 

Moreover, the maps $\gamma \mapsto P_{i,\gamma}(X)$ are coboundaries for $i>0$ and $\gamma \mapsto  P_{0,\gamma}(X)$ is a $1$-cocycle whose class in $H^1(G,V_{k-2})$ is equal to $[g]$.
\end{proposition}

 \begin{proof} By formula~\eqref{cocycle-dependence-on-z0} we write $\psi_{z_0}(\gamma)$ as a sum of $\psi_{z_1}(\gamma)$ and the coboundary $(1-\gamma)\int_{z_0}^{z_1} (X-z)^{k-2} g(z) dz$. To analyze dependence of the latter integral on $z_0$ we write it as 
\[\bal
\int_{z_0}^{z_1} (X-z)^{k-2} g(z) dz = & a_0 \int_{z_0}^{z_1} (X-z)^{k-2} dz + \int_{i \infty}^{z_1} (X-z)^{k-2} (g(z)-a_0) dz  \\
& + \int_{z_0}^{i \infty} (X-z)^{k-2} (g(z)-a_0) dz,\\
\eal\]
where $a_0$ is the constant term of the Fourier expansion of $g(z)$. In the right-hand side the first integral can be evaluated explicitly, the middle integral is independent of $z_0$ and the last summand is $o(\mathrm{Im}(z_0)^{-N})$ for all $N \ge 1$. Therefore we conclude that $\int_{z_0}^{z_1} (X-z)^{k-2} g(z) dz = \sum_{i=0}^{k-1} Q_{i,\gamma}(X) z_0^i + o(\mathrm{Im}(z_0)^{-N})$ for all $N \ge 1$ and some independent of $z_0$ polynomials $Q_{i,\gamma}\in V_{k-2}$, $0 \le i \le k-1$. Expression~\eqref{cocycle-asymptotics} now follows with 
\[
P_{0,\gamma} = \psi_{z_1}(\gamma) + (1-\gamma)Q_{0,\gamma}, \quad P_{i,\gamma} = (1-\gamma)Q_{i,\gamma}, \quad i > 0.
\]
The latter formula also proves the second claim.
\end{proof}

The above proposition shows that the coefficient in $\psi_{z_0}(\gamma)$ at every power of $X$ is regularizable as $z_0 \to \infty$  in the sense of Definition~\ref{regularized integral}. We denote 
\[
\psi_{\infty}(\gamma):= P_{0,\gamma}(X) = \sum_{i=0}^{k-2} (-1)^i\binom{k-2}{i} X^{k-2-i} \Big[ \int_{z_0} ^{\gamma z_0} z^i g(z) dz \Big]_{reg}.
\]
Since $\psi_\infty: G \to V_{k-2}$ is a 1-cocycle, we have the following relations among regularized integrals:
\begin{equation}\label{1 cocyle relation}
\psi_\infty(\gamma_1 \gamma_2) = \gamma_1 \, \psi_\infty(\gamma_2)  + \psi_\infty(\gamma_1), \qquad \gamma_1,\gamma_2 \in G.
\end{equation}

\bigskip

We now return to the computation of regularized integrals~\eqref{reg-int-to-do}. For this purpose it will be convenient to write $g(z)$ using the standard Eisenstein series 
\be\label{Eisenstein-series}
G_{ 3, \chi}(z) = \frac{1}{2} L (\chi, -2) + \sum_{ n \geq 1} \left( \sum_ { m \mid n } \chi(m) m^2 \right) q^n  \; \in M_3(\Gamma_0(N),\chi)\ee
 defined for an odd Dirichlet character $\chi: (\Z/N\Z)^\times \to \C^\times$, see~\cite[\S 4.5]{DS}. These expressions are given in Table~\ref{g-via-G}, where in all cases except {\bf D} we use the Dirichlet character $\chi$ given in Table~\ref{Z-data}. The respective L-values $L(\left(-3/\cdot\right),-2)=2/9$, $L(\left(-4/\cdot\right),-2)=-1/2$ and $L(\chi_5,-2)=-4/5-i 2/5$ were computed using MAGMA calculator, \cite{MAGMA}.

\bigskip
  
\begin{table}
\caption{ Expressions for $g(z)$ via Eisenstein series~\eqref{Eisenstein-series}.}  
\begin{center}
\begin{tabular}{c|c|l|c}\label{Eisentein expression}
& Case & $g(z)$ & b\\
\hline
&&&\\
 & \bf{A} & $ g_A(z)=1 -q -5q^2 -q^3 + 11q^4 +24q^5- \cdots$ & 6\\
 && $ \qquad\quad = -G_{3, \chi}(z) - 8 G_{3, \chi}(2z), $ & \\
 && where $\chi = \left( \frac{-3}{.}\right)$ is the non-trivial &\\
 && Dirichlet character modulo $3$ &\\
 &&&\\
 & \bf{C}  & $g_C(z) = g_A(z)$ & 6\\
 &&&\\
 &  \bf{F} & $ g_F(z)=1 +q -5q^2 +q^3 + \cdots = g_A(z+\frac12)$ & 12\\
 &&& \\
 \hline
 &&\\
  & \bf{E} & $ g_E(z)=1 -4q^2 -4q^4 + \cdots = - 4 G_{ 3, \chi}(2z) $ & 8\\
  && where $\chi = \left( \frac{-4}{.}\right)$ is the non-trivial &\\
 && Dirichlet character modulo $4$ &\\
&&&\\
& \bf{G} & $ g_G(z)= g_E(z)$ & 8 \\
  &&&\\
 \hline
 &&&\\
  & \bf{D} & $ g_D(z)= 1 - 2 q - 6 q^2 + 7 q^3 + 26 q^4 + \ldots$ & 5 \\
&& $ \qquad\quad = (-1+\frac{i}2)G_{3, \chi}(z) + (-1-\frac{i}2)G_{3, \overline{\chi}}(z),$  &\\
&& where $\chi=\chi_5$ is the odd Dirichlet character &\\
&& modulo $5$ with $\chi_5(2)=i$ & \\  
\end{tabular}
\end{center}
\label{g-via-G}  
\end{table}

\bigskip
 
We denote 
$$
\tilde{g}(z) = \sum_{ n \geq 0}  \tilde{a}_n q^n := g(z+ 1/ b) = 1 + \sum_{n \geq 1} a_n   \exp(2 \pi i / b)^n q^n.   
$$ 
Note that from the above representation in terms of Eisenstein series it follows that $|\tilde{a}_n|=|a_n| = O(n^2)$ as $n \to \infty$ and therefore the  twisted Dirichlet series $L(\tilde{g},s) = \sum_{n \ge 1}  \frac{\tilde{a}_n}{n^s}$ is convergent for $\mathrm{Re}(s)>3$. Below we will evaluate it explicitly    in some of our cases, and from these expressions it will follow  that   $L(\tilde{g},s)$ is holomorphic for $s \in \C$ except for $s=3$ where it has a simple pole.

\begin{lemma} \label{inverse derivative} Suppose that $L(\tilde g, s)$ has a meromorphic continuation to $\C$ with a possible pole only at $s=3$. As $ |z_0| \to \infty$ one has
\begin{equation}\label{reg-int-via-L}
  \int_{ z_0} ^{\frac{z_0}{bz_0+1}} g(z) dz = \frac{1}{b}-z_0+ \frac{L(\tilde{g}, 1)}{2 \pi i}- \frac{b^2 ( b z_0 +1)^2}{8 \pi^3 i}\mathrm{Res}_{s=3} L (\tilde{g}, s) +o(1).
  \end{equation}
\end{lemma}  
  
\begin{proof} 
We will do the change of variables $z=w+\frac1b$. Denote $\tilde g^*(z)=\tilde g(z) - a_0 = \tilde g(z) - 1$. Put $w_0= z_0 -\frac{1}{b}$. We then have 
\begin{align} \label{integral_g(z)}
   \int_{z_0}^{\frac{z_0}{bz_0+1}} g(z) dz & = \int_{w_0}^{-\frac{1}{b(b w_0+2)}} \tilde{g}(w) dw= -\frac{1}{b( b w_0+2)}- w_0  + \int_{w_0}^{-\frac{1}{b(b w_0+2)}} \tilde{g}^*(w) dw {\color{red}  \nonumber }\\
   &= -w_0 + \frac{1}{2 \pi i} F \left(-\frac{1}{b(bw_0+2)}\right) +o(1), \mbox{ for }  |w_0| \rightarrow  \infty, 
\end{align}
where
\[
F(z) = \sum_{n \ge 1} \frac{\tilde a_n}{n}q^n
\] 
is the function satisfying $\frac1{2\pi i} \frac{d}{dz} F = q \frac{d}{dq} F = {\tilde g}^*$.  We will now show that 
\begin{equation}\label{antiderivative-asymptotics}
 F(it) = L (\tilde{g}, 1) + \frac{1}{(2 \pi t)^2 } \mathrm{Res}_{s=3} L (\tilde{g}, s) +o(1), \; t \to 0_+.
\end{equation}
Formula~\eqref{reg-int-via-L} then follows from this and~\eqref{integral_g(z)}. Recall that the Dirichlet series for $L(\tilde g, s)$ is convergent for $\mathrm{Re}(s)>3$. Therefore for $\mathrm{Re}(s)>2$ we have
\[
\frac{\Gamma(s)}{(2\pi)^s} L(\tilde g, 1 + s) = \frac{\Gamma(s)}{(2\pi)^s} \sum_{n \ge 1}\frac{\tilde a_n}{n^{s+1}} = \sum_{n \ge 1}\frac{\tilde a_n}{n} \int_0^\infty t^{s-1} e^{-2 \pi n t} dt = \int_0^{\infty} t^{s-1} F(i t) dt.
\]  
Since $|\Gamma(s)|$ is small when $|\mathrm{Im}(s)|$ is large and $L(\tilde g, 1+s)$ is uniformly bounded when $\mathrm{Re}(s) \ge 2 + \varepsilon$ with any fixed $\varepsilon > 0$, we can apply the Mellin inversion theorem. One recovers
\[
F(i t) = \frac1{2 \pi i} \int_{c-i \infty}^{c+i \infty} \frac{\Gamma(s)}{(2\pi t)^s} L(\tilde g, 1 + s) ds
\] 
with any real $c > 2$. If we integrated over the vertical line with $c<0$, then the integral would be $o(1)$ when $t \to 0_+$. This is because $|t^{-s}|=t^{- \mathrm{Re}(s)}$. We can move the line of integration from $c>2$ to $-1<c<0$. The integrand has two poles between those lines: $s=0$ is a pole for $\Gamma(s)$ and $s=2$ is a pole for $L(\tilde g,1+s)$. Therefore one needs to add the residues of the integrand at those points, which yields formula~\eqref{antiderivative-asymptotics}.
\end{proof}    

\begin{lemma}\label{case A, C} In case {\bf A} one has

\be\label{L-A}\bal
L(\tilde{g}, s) = &-\frac{1}{2}  (1-3^{1-s})(1-2^{2-s})(1+2^{1-s}) \zeta(s) L( \chi, s-2) \\
& - i \frac{\sqrt{3}}{2} (1-3^{2-s})(1-2^{1-s})(1+2^{2-s})  \zeta(s-2) L(\chi, s).
\eal\ee
 For this  twisted Dirichlet series we have 
\[
L (\tilde{g}, 1) = 0 \quad\text{  and }\quad \mathrm{Res}_{s=3} L(\tilde{g}, 3) = -i\frac{ \pi^3}{54}.
\]
\end{lemma}

{
\begin{proof} From the expression for $g(x)$ in Table~\ref{g-via-G} we find that its Fourier coefficients are given by $a_n = \sum_{d|n} (-1)^d \chi(d) d^2$. Denote $\alpha=\exp(2 \pi i /6)$,
\[
A_k(s) = \sum_{m \geq 1} \frac{\alpha^{k m}}{m^s} \qquad{  and } \qquad B_k(s) = \sum_{m \geq 0} \frac{1}{(6 m + k)^s}
\]
for $k \ge 1$. For the twisted  Dirichlet series we then have the following 
\be\label{twisted-L}\bal
L(\tilde{g}, s) &= \sum_{d \ge 1}\sum_{m \ge 1} \frac{(-1)^d \chi(d) d^2 \alpha^{md}}{(md)^s} = \sum_{d \ge 1} \frac{(-1)^d \chi(d)}{d^{s-2}} A_d(s) \\
&= - B_1(s-2)A_1(s) - B_2(s-2)A_2(s) + B_4(s-2)A_4(s) + B_5(s-2)A_5(s). 
\eal\ee  
We now list some observations about the $A$'s:
\[\bal
&A_1(s)+A_4(s) = \sum_{m=1}^\infty \frac{\alpha^m}{m^s} + \sum_{m=1}^\infty (-1)^m \frac{\alpha^m}{m^s} = 2 \sum_{2 | m} \frac{\alpha^m}{m^s} = 2^{1-s} A_2(s),  \\
&A_2(s)+A_5(s) = \sum_{m=1}^\infty \frac{\alpha^{2m}}{m^s} + \sum_{m=1}^\infty (-1)^m \frac{\alpha^{2m}}{m^s} = 2 \sum_{2 | m} \frac{\alpha^{2m}}{m^s} = 2^{1-s} A_4(s).\\
\eal\] 
In~\eqref{twisted-L} we will substitute $A_1(s)=2^{1-s}A_2(s)-A_4(s)$ and $A_5(s)=2^{1-s}A_4(s)-A_2(s)$. Since $\alpha^2 = -\frac12+i \frac{\sqrt{3}}2$, $\alpha^4 = -\frac12- i \frac{\sqrt{3}}2$ and $\alpha^6=1$, we obtain that
\[\bal
A_2(s) &= -\frac12 \sum_{3 \nmid m}\frac1{m^s} + i \frac{\sqrt{3}}2 L(\chi,s) + \sum_{3 | m}\frac1{m^s} \\
&= -\frac12 (1-3^{-s})\zeta(s) + i \frac{\sqrt{3}}2 L(\chi,s) + 3^{-s}\zeta(s) \\
&= -\frac12 (1-3^{1-s})\zeta(s) + i \frac{\sqrt{3}}2 L(\chi,s)\\
\eal\]
and similarly
\[
A_4(s) = -\frac12 (1-3^{1-s})\zeta(s) - i \frac{\sqrt{3}}2 L(\chi,s).
\]
Substituting these expressions into~\eqref{twisted-L} we obtain
\be\label{L-A-intermediate}\bal
L(\tilde{g}, s) &= - B_1(s-2)\left(2^{1-s}A_2(s)-A_4(s)\right) - B_2(s-2)A_2(s) \\
& \quad + B_4(s-2)A_4(s) + B_5(s-2)\left(2^{1-s}A_4(s)-A_2(s)\right)\\
&= \left( -2^{-1-s} B_1 - B_2-B_5 \right)(s-2) \, A_2(s) \\
&\quad + \left( B_1 + B_4+ 2^{-1-s}B_5 \right)(s-2)\, A_4(s) \\
&= \left( -(1-2^{-1-s}) (B_1 - B_5) + (B_2-B_4)\right)(s-2) \times \frac12 (1-3^{1-s}) \zeta(s) \\
&= \left( -(1+2^{-1-s}) (B_1 + B_5) - (B_2+B_4)\right)(s-2) \times i\frac{\sqrt{3}}2 L(\chi,s).\\
\eal\ee
It remains to do a few observations about $B$'s. One can immediately see that
\[
B_2(s)-B_4(s) = 2^{1-s} L(\chi,s), \qquad B_1(s)-B_5(s) = (1+2^{-s}) L(\chi,s).
\]
Substituting these formulas into the first row of our latest expression for $L(\tilde{g},s)$ yields the first row in formula~\eqref{L-A}. Now observe that $B_1(s)+B_3(s)+B_5(s)=\sum_{m \ge 1, 2 \nmid m} m^{-s}=(1-2^{-s})\zeta(s)$ and $B_3(s)=3^{-s}\sum_{m \ge 1, 2 \nmid m} m^{-s} = 3^{-s}(1-2^{-s})\zeta(s)$. This implies
\be\label{B1+B5}
B_1(s)+B_5(s) = (1-2^{-s})(1-3^{-s})\zeta(s).
\ee
Subtracting this from $B_1(s)+B_2(s)+B_4(s)+B_5(s)=\sum_{m \ge 1, 3 \nmid m} m^{-s}=(1-3^{-s})\zeta(s)$ we obtain 
\be\label{B2+B4}
B_2(s)+B_4(s) = 2^{-s}(1-3^{-s})\zeta(s).
\ee
Substituting~\eqref{B1+B5} and~\eqref{B2+B4} into the last row in~\eqref{L-A-intermediate} we then obtain the second row in formula~\eqref{L-A}. This completes the proof of~\eqref{L-A}.

Since $L(\chi,s)$ is entire and $\zeta(s)$ is holomorphic everywhere except for a simple pole at $s=1$, we find that the first row in~\eqref{L-A} is entire and the second row has a unique simple pole at $s=3$. To evaluate the Dirichlet series $ L(\tilde{g},s)$ at $s=1$ we note that $L(\chi,-1)=0$, see~\cite[ Theorem 4.2]{WaCF}. Therefore  there is no contribution from the first row of~\eqref{L-A}, and the term $(1-2^{1-s})$ in the second row causes it to vanish at $s=1$. Therefore $L(\tilde{g},1)=0$. Since $\mathrm{Res}_{s=1}\zeta(s)=1$, we compute that
\[\bal
\mathrm{Res}_{s=3} L(\tilde{g},s) &= - i \frac{\sqrt{3}}{2} (1-3^{-1})(1-2^{-2})(1+2^{-1}) L(\chi, 3)\\
&= - i \frac{\sqrt{3}}{2}  \times \frac23 \times \frac34 \times \frac32 \times \frac{4 \pi^3 \sqrt{3}}{243} = - i \frac{\pi^3}{54}.
\eal\]   
Here we substituted the value $L(\chi,3)=\frac{4 \pi^3 \sqrt{3}}{243}$ which follows from $L(\chi,-2) = L(\left(-3/ \cdot \right), -2) =-\frac29$ and the functional equation of $L(\chi,s)$ with respect to $s \leftrightarrow 1-s$, see~\cite[ Chapter 4]{WaCF}.     
\end{proof}
}

\medskip

\begin{lemma} \label{case D} For case {\bf D}, with character $\chi$ given in Table~\ref{Eisentein expression}, one has

\be\label{statementD}\bal L(\tilde{g}, s) = 
&   \; -L(2 \mathrm{Re}(\chi) + \mathrm{Im}(\chi), s-2) \, 5^{-s} \zeta(s) \\
&+\frac12 \cos\left(\frac{2\pi}5\right) \Big[- L(2 \mathrm{Re}(\chi)+\mathrm{Im}(\chi), s-2) \zeta(s)( 1 - 5^{-s}) \\
& \qquad\qquad\quad  - L(2 \mathrm{Re}(\chi)-\mathrm{Im}(\chi), s-2) \, L( \chi^2, s)  \Big]  \\
&+\frac12 \cos\left(\frac{4\pi}5\right) \Big[- L(2 \mathrm{Re}(\chi)+\mathrm{Im}(\chi), s-2) \zeta(s)( 1 - 5^{-s}) \\
& \qquad\qquad\quad  + L(2 \mathrm{Re}(\chi)-\mathrm{Im}(\chi), s-2) \, L( \chi^2, s)  \Big]  \\
&+\frac{i}2 \sin\left(\frac{2\pi}5\right) \Big[- L(2 \mathrm{Re}(\chi)-\mathrm{Im}(\chi), s) \, \zeta(s-2)( 1 - 5^{2-s}) \\
& \qquad\qquad\quad  - L(2 \mathrm{Re}(\chi)+\mathrm{Im}(\chi), s) \, L( \chi^2, s-2)  \Big]  \\
&+\frac{i}2 \sin\left(\frac{4\pi}5\right)\Big[- L(\mathrm{Re}(\chi)+ 2 \mathrm{Im}(\chi), s) \, \zeta(s-2)( 1 - 5^{2-s}) \\
& \qquad\qquad\quad  + L(\mathrm{Re}(\chi)- 2 \mathrm{Im}(\chi), s) \, L( \chi^2, s-2)  \Big] . \\
\eal  \ee
This Dirichlet series is holomorphic for all $s \in \C$ except for $s=3$, where it has a simple pole. 
 We have
\begin{align*}
L (\tilde{g}, 1) 
%& = i \sin \left(\frac{ 2 \pi}{5}\right)\left[ \frac{Re(L(\chi,1)}{15} + \frac{11 Im(L(\chi, 1))}{30} \right]+ i \sin \left( \frac{4 \pi}{5}\right)\left[ \frac{Im(L(\chi,1)}{15} - \frac{11 Re(L(\chi, 1))}{30} \right]\\
& =  -\frac{\pi}{60}i \quad\mbox{ and }\quad 
 \mathrm{Res}_{s =3} L (\tilde{g}, s)
 %& = \frac{4i}{5} \sin\left(\frac{ 2 \pi}{5}\right) \left[-Re(L(\chi,3) + \frac{ Im(L(\chi, 3))}{2} \right]-\frac{4i}{5} \sin\left( \frac{4 \pi}{5}\right) \left[ Im(L(\chi,3) + \frac {Re(L(\chi, 3))}{2} \right]\\
  = -\frac{4\pi^3}{125}i.
 \end{align*}
 \end{lemma}

\begin{proof}

From the expression for $g(z)$ in Table~\ref{g-via-G} we find that its Fourier coefficients are given by $a_n = \sum_{d|n} \tau(d) d^2$, where $\tau = - 2 \mathrm{Re}(\chi) - \mathrm{Im}(\chi)$. Let $\alpha=\exp(2 \pi i /5)$. Denote 
\[
A_k(s) = \sum_{m \geq 1} \frac{\alpha^{k m}}{m^s} \qquad{  and } \qquad B_k(s) = \sum_{m \geq 0} \frac{1}{(5 m + k)^s}
\]
for $k \ge 1$. 
For the twisted Dirichlet series, we then have
\be\label{twisted L fun D}\bal
L(\tilde{g}, s)  & = \sum_{ d \geq 1}\sum_{ m \geq 1}  \frac{ \tau(d) d^2 \alpha^{md}}{(md)^{s}}   = \sum_{d \geq 1}\frac{ \tau(d) }{d^{s-2}} A_d(s)\\
& =  -2 B_1(s-2) A_1(s) + 2 B_4(s-2) A_4(s)-B_2(s-2) A_2(s) + B_3(s-2) A_3(s).
\eal\ee
 We express each of $B_1(s),B_2(s),B_3(s)$ and $B_4(s)$ can be expressed as a linear combination of $L(\chi^j,s)$
for $j=1,2,3,4$:
\[\bal
&B_1(s) = \frac14 L(\chi,s) + \frac14 L(\chi^2,s) + \frac14 L(\chi^3,s) + \frac14 L(\chi^4,s), \\   
&B_2(s) = -\frac{i}4 L(\chi,s) - \frac14 L(\chi^2,s) + \frac{i}4 L(\chi^3,s) + \frac14 L(\chi^4,s), \\   
&B_3(s) = \frac{i}4 L(\chi,s) - \frac14 L(\chi^2,s) - \frac{i}4 L(\chi^3,s) + \frac14 L(\chi^4,s), \\   
&B_4(s) = -\frac14 L(\chi,s) + \frac14 L(\chi^2,s) - \frac14 L(\chi^3,s) + \frac14 L(\chi^4,s). \\   
\eal\]
Note that $\chi^3=\overline{\chi}$ and $L(\chi^4,s)=(1-5^{-s})\zeta(s)$. Functions $A_k(s)$ can be also written as such linear combinations using the above formulas for $B$'s and the fact that $A_k = \sum_{j=1}^4 \alpha^{kj}B_j + 5^{-s}\zeta(s)$. With the help of a computer, we substituted all these expressions into~\eqref{twisted L fun D} and found that the result equals~\eqref{statementD}. 

For $k=1,\ldots,4$ functions $5^{s} B_k(s) = \sum_{m \ge 0}\frac1{(m+k/5)^s}$ are Hurwitz zeta functions. They are known to have analytic continuation to $\C\setminus \{ 1 \}$ and a simple pole of residue $1$ at $s=1$. Therefore $L({\rm Re}\chi,s)=B_1(s)-B_4(s)$, $L({\rm Im}\chi,s)=B_2(s)-B_3(s)$ and $L(\chi^2,s)=B_1(s)-B_2(s)-B_3(s)+B_4(s)$ are entire functions in the complex plane.  Hence the upper three rows of expression~\eqref{statementD} may only have a pole at $s=1$ and the lower two rows may only have a pole at $s=3$. Using MAGMA we check that $L(\chi,-1)=L(\overline\chi,-1)=0$, from which it follows that 
\[
L({\rm Re} \chi, -1) = L({\rm Im} \chi, -1) = 0. 
\]
By this reason the upper three rows of~\eqref{statementD} are entire functions. To compute the residue of the lower two rows at $s=3$ we need the values $L({\rm Re} \chi,3)= {\rm Re}L(\chi, 3)$ and $L({\rm Im} \chi,3)={\rm Im} L(\chi, 3)$. The value $L(\chi, 3)$ can be computed from the functional equation satisfied by $L(\chi,s)$ using the value $L(\chi,-2)=-4/5-2/5 i$, which we found earlier using MAGMA. Thus we find that
\[\bal
&L({\rm Re} \chi,3) = \frac{8 \pi^3}{625}\left(2 \sin\left(\frac{2\pi}5\right) + \sin\left(\frac{4\pi}5\right)\right),\\
&L({\rm Im} \chi,3) = \frac{8 \pi^3}{625}\left(- \sin\left(\frac{2\pi}5\right) + 2 \sin\left(\frac{4\pi}5\right)\right).\\
\eal\]    
The residue of $\zeta(s-2)$ at $s=3$ equals $1$, and therefore we compute the residue of $L(\tilde g,s)$ at $s=3$ as the sum of values of remaining terms in the last two rows of our expression, with $\zeta(s-2)$ removed from them:
\[\bal
\mathrm{Res}_{s =3} L (\tilde{g}, s) = \frac{i}{2} \sin\left(\frac{2\pi}5\right) \left(-2 L({\rm Re} \chi,3) +  L({\rm Im} \chi,3)\right) \cdot\left(1-\frac15\right)\\
+ \frac{i}{2} \sin\left(\frac{4\pi}5\right) \left(-L({\rm Re} \chi,3) - 2  L({\rm Im} \chi,3)\right) \cdot\left(1-\frac15\right) \\
= \frac{i}{2} \cdot 
\frac{8 \pi^3}{625}\left(-5 \sin\left(\frac{2\pi}5\right)^2 - 5 \sin\left(\frac{4\pi}5\right)^2\right)\cdot \frac 45 = - \frac{4 \pi^3}{125} i.
\eal\]
Here we used the identity $\sin(2\pi/5)^2+\sin(4\pi/5)^2=5/4$, which can be proved using elementary transformation rules for trigonometric functions.

It remains to determine $L(\tilde g,1)$. We first show that the contribution from the first three rows of~\eqref{statementD} vanishes. Since $L({\rm Re} \chi,-1)=L({\rm Im} \chi, -1)=0$, the summands with $L(\chi^2,s)$ do not contribute to the sum at $s=1$. Therefore the value of the sum of the first three rows equals to the limit $-L(2 {\rm Re}(\chi)+{\rm Im}(\chi),s-2)\zeta(s)$ as $s\to 1$ times
\[
\frac15+\frac12\cos\left(\frac{2\pi}{5}\right)\left(1-\frac15\right)+\frac12\cos\left(\frac{4\pi}{5}\right)\left(1-\frac15\right)=0.
\] 

To compute the value of the sum of the last two rows in~~\eqref{statementD} at $s=1$, we use the values $\zeta(-1)=-1/12$, $L(\chi^2, -1) = -2/5$ and
 $  L(\chi, 1) = \left(\frac{2\pi}{25} + \frac{6\pi}{25}i \right)\sin\left( \frac{4\pi}{5}\right)+\left(\frac{6\pi}{25}- \frac{2\pi}{25}i \right)\sin\left(\frac{2\pi}{5}\right)$. After putting all these values into our expression, we get the desired value of the Dirichlet series $L(\tilde{g},s)$ at $s=1$.
\end{proof}

\noindent{\it Proof of Theorem \ref{kappa2-thm}.}
\smallskip

(Case {\bf A} and {\bf C}:) %As above, we denote $ g_A(z) = g_C(z) =: g(z)$ in this case. 
By  Proposition \ref{case A, C}, we have the values of $ L(\tilde{g}, 1)$ and $ \mathrm{Res}_{s=3} L(\tilde{g}, 3)$. Putting these values 
 in Proposition \ref{inverse derivative}, we get 
\be\label{asymptotics-A}
 \int_{z_0} ^{\frac{z_0}{6 z_0+1}} g(z) dz =  \frac{1}{4}+ 3 z_0^2 + o(1) \mbox{ as } z_0  \rightarrow i \infty.
 \ee Therefore, we conclude $\displaystyle \left[ \int_{z_0} ^{ \frac{z_0}{6z_0 +1}} g(z) dz \right] _{reg} = \frac{1}{4}$ and by Theorem~\ref{kappa-1-is-0-thm} the value of $\kappa_2$ equals $-\frac{\pi^2}6$.
 
 \bigskip

 (Case {\bf D}:) Following a similar line of argument as above, in this case, 
 $$\left[ \int_{z_0} ^{ \frac{z_0}{5z_0 +1}} g_D(z) dz \right] _{reg} = \frac{7}{24}$$ and by Theorem~\ref{kappa-1-is-0-thm} the value of $\kappa_2$ equals $-\frac{7}{30}\pi^2$.
 
 \bigskip

 (Case {\bf F}:) Performing the change of variable $z=u-\frac12$ we obtain 
 \be\label{caseF}
 \int_{z_0}^\frac{z_0}{12z_0 +1} g_F(z) dz = \int_{u_0}^{ \frac{7 u_0 -3}{12 u_0 -5}} g(z) dz,
 \ee
 where $g(z)=g_{A}(z)$ is the modular form on $\Gamma_1(6)$ as in the previous considered case. We will use Proposition~\ref{regularization-prop} to find the asymptotics of this integral as $u_0 \to \infty$. Let $G=\Gamma_1(6)$. One can check that this group is generated by the parabolic elements   
\[ 
\sigma_0 = \begin{pmatrix}1 & 1 \\ 0 & 1\end{pmatrix}, \qquad \sigma_c =  \begin{pmatrix}1 & 0 \\ 6 & 1\end{pmatrix},
\qquad \sigma_{c'} = \begin{pmatrix}
7 & -3 \\
12 & -5
\end{pmatrix}
\]
fixing the cusps $\lambda = \infty, 0$ and $\frac12$ respectively. Following the notation of Proposition~\ref{regularization-prop} for $\gamma \in G$ we write
\[
\psi_{z_0}(\gamma) = \int_{z_0}^{\gamma z_0} (X-z)g(z) dz = \sum_{i=0}^2 P_{i,\gamma}(X) z_0^i + o(1), \quad z_0 \to \infty.
\]
Since
\[
\psi_{z_0}(\sigma_0) = \int_{z_0}^{z_0+1} (X-z) g(z) dz = X - \frac{(z_0+1)^2-z_0^2}{2} + o(1), \quad z_0 \to \infty
\]
we have
\be\label{sigma-0-values}
P_{0,\sigma_0}(X) = X-\frac12, \quad P_{1,\sigma_0}(X) = -1, \quad P_{2,\sigma_0}(X) = 0.
\ee
We note that if $\sigma \in G$ is a parabolic element fixing a cusp $\lambda \in \Q$ and the form $(z-\lambda)g(z) dz$ is bounded near $z = \lambda$ then for any $z_0 \in \sH$ one has
\[
\int_{z_0}^{\sigma z_0} (z - \lambda) g(z) dz = 0. 
\]
Indeed, this integral is independent of $z_0$ because the form is $\sigma$-invariant and we send $z_0 \to \lambda$ to show that the integral vanishes. Therefore in the discussed case one has 
\be\label{asymptotics-parabolic}
\psi_{z_0}(\sigma) = \int_{z_0}^{\sigma z_0} (X-z)g(z) dz = (X - \lambda)\int_{z_0}^{\sigma z_0} g(z) dz.\\
\ee
Boundedness of $(z-\lambda)g(z) dz$ near $z=\lambda$ depends only on the $G$-orbit of $\lambda$. There are 4 $G$-orbits of cusps for our group, they are represented by $\lambda=\infty, 0, \frac12$ and $\frac13$, and we have boundedness for the latter three orbits. Formulas~\eqref{asymptotics-A} and~\eqref{asymptotics-parabolic} for $\sigma=\sigma_c$ yield  
\be\label{sigma-c-values}
P_{0,\sigma_c}(X) = \frac14 X, \quad P_{1,\sigma_c}(X) = 0, \quad P_{2,\sigma_c}(X) = 3 X.
\ee
Recall that for $i = 1, 2$ the maps $ \gamma \mapsto P_{i, \gamma}(X)$ are $1$-coboundaries, that is there exist polynomials $Q_1(X), Q_2(X)$ of degree $\le 1$ such that $P_{i, \gamma} = (1-\gamma) Q_i$ for $i=1,2$. From~\eqref{sigma-0-values} and~\eqref{sigma-c-values} we recover that $Q_1(X) = -X$ and $ Q_2(X) = \frac{1}{2}$. With this we obtain $P_{1,\sigma_{c'}}(X)=(1-\sigma_{c'})Q_1 = -6 X + 3$ and $P_{2,\sigma_{c'}}(X)=(1-\sigma_{c'})Q_2 = 6 X - 3$. In the view of~\eqref{asymptotics-parabolic} for $\sigma=\sigma_{c'}$ we have 
\be\label{sigma-c'-values}
P_{0,\sigma_{c'}}(X) = \beta(X-\frac12), \quad P_{1,\sigma_{c'}}(X) = -6(X-\frac12), \quad P_{2,\sigma_{c'}}(X) = 6 (X-\frac12),
\ee
with the unknown constant $\beta = \left[ \int_{z_0}^{\sigma_{c'}z_0} g(z) dz \right]_{reg}$. We will now determine $\beta$ using the $1$-cocycle property of $\gamma \mapsto P_{0, \gamma}(X)=\psi_{\infty}(\gamma)$. To this end we note that $\sigma=\sigma_{c} \sigma_0^{-1} \sigma_{c'} = \begin{pmatrix}-5 & 2 \\ -18 & 7\end{pmatrix}$ is a parabolic element fixing the cusp $\lambda = \frac13$. Using the cocycle relation~\eqref{1 cocyle relation} we obtain
\[\bal
\psi_{\infty}(\sigma) &= \psi_{\infty}(\sigma_{c}) + \sigma_{c} \psi_\infty(\sigma_0^{-1}) + \sigma_{c} \sigma_0^{-1} \psi_\infty(\sigma_{c'})  \\
&= \psi_{\infty}(\sigma_{c}) + \left(\psi_\infty(\sigma_{c'}) - \psi_\infty(\sigma_0)\right)| \sigma_0 \sigma_{c}^{-1} \\
&= \frac{X}{4} + (\beta-1) ( X-\frac12) \Big| \begin{pmatrix} -5 & 1 \\ -6 & 1\end{pmatrix} = (-2 \beta + \frac94)X + \frac12(\beta - 1).
\eal\]
In the view of~\eqref{asymptotics-parabolic} this polynomial should be a multiple of $(X-\frac13)$, from which we find that $\beta=\frac32$. With this value of $\beta$, formula~\eqref{sigma-c'-values} yields
\[
\int_{z_0}^{\sigma_{c'}z_0} g(z) dz = \frac32 - 6 z_0 + 6 z_0^2 + o(1), \quad z_0 \to \infty.  
\]
Using this asymptotics in the right-hand side of~\eqref{caseF} we find that $\displaystyle \left[ \int_{z_0} ^{ \frac{z_0}{12z_0 +1}} g_F(z) dz \right] _{reg} = 0$. It follows that $\kappa_2=0$.

 \bigskip
 
 (Cases {\bf E} and {\bf G}:) Since $g_E=g_G$ and $b=8$ in both cases, the values of $\kappa_2$ should be also equal. For the case {\bf G} we already computed in~\eqref{kappa values of G} that $\kappa_2=-\frac{\pi^2}{12}$.
 \qed

\end{document}